\documentclass{notices}
\usepackage{amsfonts,amssymb,amsmath,amscd}
\usepackage{tikz}
\usetikzlibrary{trees}

\title{
A Field Guide to Ethics in Mathematics
}

\author{
  Maurice Chiodo
  \affil{ 
    Maurice Chiodo is a research associate at the Centre for the Study of Existential Risk, University of Cambridge, and the principal investigator of the Ethics in Mathematics Project. His email address is mcc56@cam.ac.uk
    }
  \and
  Dennis M\"uller
  \affil{
    Dennis M\"uller is a graduate student at RWTH Aachen University, Germany, and a research affiliate at the Centre for the Study of Existential Risk, University of Cambridge. His email address is dennis.mueller3@rwth-aachen.de
   }
}

\begin{document}

\maketitle

\newpage
\section*{Introduction}
Mathematics has become inescapable in modern, digitized societies: there is hardly any area of life left that isn't affected by it, and we as mathematicians play a central role in this. Our actions affect what others, in particular our students, decide to do with mathematics, and how mathematics affects the world\cite[p.~2]{porter2023nonexperts}, for better or worse. In return, the study of ethics in mathematics (EiM) has become increasingly important, even though it is still unknown to many. This exposition tries to change that, by motivating ethics in mathematics as an interesting, tractable, non-trivial, well-defined and good research area for mathematicians to consider. 

\subsection*{What is ethics in mathematics?}
Historically, many of those working on ethics in mathematics have focused on questions related to specific mathematical subdisciplines (e.g.~finance, cryptography or statistics), but over time the field has begun to ask questions that are fundamental to any area of mathematics\cite{Muller.2022}. Today, this emerging area of research studies the moral principles of all mathematical practice, addressing the ethical questions related to its applications, teaching, and research, by focusing on the responsibilities of those engaged in mathematics. Ethics in mathematics usually builds on the assumption that the majority of ``mathematicians want to be ethical even when [some of them] don’t think ethics apply to them''\cite[p.~4]{buell2022leveraging}.

This work often transcends the boundaries of mathematics, and involves consulting not just the obviously relevant disciplines such as philosophy and law, but also the political and social sciences and psychology\cite{chiodo2023manifesto}. Consider, for example, recent advances in automated theorem proving by artificial intelligence. ``What constitutes a proof?'' and ``Who counts as an author?'' aren't merely practical questions; they're also ethical concerns with deep psychological and social components: Are we prepared to hand over authorship to an AI \textit{and trust it}? Similarly, in applied mathematics, the use of approximative formulas also leads to complicated ethical issues. For example, the use of Gaussian copulas to model dependencies between random variables  played a fundamental role in the Global Financial Crisis of 2008, leading to some authors calling it the ``formula that killed Wall Street''\cite{MacKenzie.2014}. Studying the ethics of using approximate formulas to model complex situations requires us to look beyond mathematics, and use methodologies from many other disciplines.

Mathematics affects our modern societies and our planet with an ever-increasing impact. Today, we cannot hide from mathematics. Someone will have done mathematics that affected us in one way or another, and our mathematics will have affected someone else - positively or negatively. Ethics in mathematics studies this, and guides mathematicians through this jungle of problems. It is, ultimately, an area of philosophy of mathematics, that deeply values the ``problems emerging from the actual mathematical practice'' and its people\cite[p.~115]{Skovsmose.2023}.

\subsection*{Ethics in mathematics and other efforts} The difference between ethics in mathematics and other similar endeavors appears to be mostly a question of perspective and focus. 
For example, many researchers working on ``Mathematics for Social Justice'' are predominantly focused on (solving) problems surrounding power and identity, opportunities and injustices, including economic, judicial, racial and other privileges, in relation to mathematics and its teaching. For example, they might study the role that mathematics plays in (countering) gerrymandering, or how we can create better classroom environments for teaching. On the other hand, while ethics in mathematics is certainly concerned with these questions, it also studies ethical questions not restricted to balance, equity and fairness. Over time, the different foci have led to complementary methodologies in research and teaching\cite{müller2022situating}, with many fruitful interactions between them. Questions of social justice are often related to other ethical questions, and almost all ethical questions have a social justice component. All of these can be understood within the context of a sociopolitical turn in mathematics and its education\cite{Gutierrez.2013}, claiming that mathematics is not a neutral and pure player in (modern) societies\cite{Ernest2020}.

\subsection*{The structure of ethics in mathematics}

Embedded throughout this paper are  examples from three (partially overlapping) areas of concern:

\begin{itemize}
\setlength\itemsep{-1pt}
    \item the ethics of ensuring that mathematics remains a stable and continuing body of knowledge,
    \item the ethical and social questions surrounding our mathematical, and other, communities, and 
    \item the impact of mathematics on wider society and our planet. 
\end{itemize}
Most questions in ethics in mathematics fall into one or more of these three areas, and we can use these as a starting point for thinking about the ethics of our, and others', mathematical work and research. Does our mathematical work make mathematics a better area of research? Is the way we do and teach mathematics good for the mathematical community? Does our work have a positive impact on wider society?  For industrial applications of mathematics, the recent ``Manifesto for responsible development''\cite{chiodo2023manifesto} acts like a handbook to provide a structured list guiding the reader through the most common ethical questions arising from the use of mathematics. 

But, as shown later in this article, these areas of concern are not restricted to applied mathematics. Pure mathematics can give rise to ethical concerns, which include: Is the motivation that goes into a piece of mathematics good? Are we writing down our work in a way that other researchers can understand, or are we trying to fudge the boundaries because we want to hide problems? Is the way we do and teach mathematics inclusive to others, enhancing their understanding? Are we working on problems that matter (to us or someone else)? How can our (pure) mathematics be used? And are there areas of dual-use we should be concerned about? Regardless of the type of mathematics we do, be it pure or applied, it's often necessary to look beyond our area of concern to fully understand the ethics of our work.

\subsection*{How to use this article}
This article is intended to give mathematicians a brief introduction to ethics in mathematics. It's not a philosophical treatise of the subject, but rather an explanation of why doing ethics in mathematics is not so different from doing mathematics, and how we can contribute to the effort or even just use existing resources in our work and teaching. This article is a research exposition of a slightly different kind, explaining why ethics in mathematics satisfies the characteristics of a problem worthwhile of a mathematician's attention: it is interesting, tractable, non-trivial, well-defined and (morally) good. Each of these aspects is presented in one page, with the hope of making this article usable in many settings. One can always find five minutes in a busy mathematics course to talk about mathematics and its community at large, so each of the following five pages is written with the intent of being a useful handout for such a five-minute excursion.

\newpage
\section*{Why is EiM interesting?}

Many of us do mathematics because we find it interesting, and often we choose problems to work on that generate interesting mathematics. Andrew Wiles already argued that ``it's fine to work on any problem, so long as it generates interesting mathematics along the way’’ \cite{PBSNOVA2000}. At face value, ethics in mathematics might not seem like something that would pique our interest. After all, how could studying the \textit{ethics} of mathematical work generate new problems of \textit{mathematical} interest? But actually, it does.

As mathematicians, we know that theorems often create more open problems. The more we know, the more questions we can and want to pose, and the same is true for ethics in mathematics. Simply said, mathematics creates ethical questions, and in return ethical questions generate interesting mathematics. Consider, for example, the vast amount of technical research on algorithmic fairness. As algorithms are now increasingly influencing many facets of our lives, new ethical questions have emerged. But to answer these questions, we need new mathematics. ``What is fairness?’’ is a deep philosophical question, but it’s also an inherently mathematical one. It has led mathematicians and computer scientists to develop new metrics to measure individual and group fairness, connecting their area of mathematics with ethics, law and philosophy.

During this process, mathematics often helps to determine the limits of these ethical questions. When generative AI (e.g.~ChatGPT) ``hallucinates’’ and tells us something untrue, we may ask if it’s always possible to identify these untrue statements, or even better, create stable neural networks that accurately output true statements. Mathematics shows this is sometimes provably impossible\cite{Colbrook.2022}. Creative, deep and new mathematics goes hand in hand with difficult and urgent ethical questions.

But ethics in mathematics is more than the technical work of individuals. Due to its inherent interdisciplinary and often transdisciplinary nature, this young research area currently primarily consists of several small, widely dispersed communities, straddling many areas beyond the boundary of mathematics. Onboarding people for this work is an interesting and necessary challenge in itself. To properly address many of the ethical questions raised by mathematical work, we will need to come together as a community. Some questions (e.g.~ethical questions about proof standards) necessitate  community consensus, but all questions benefit from productively assembling our unique knowledge, experience, and perspective. Bringing these separate communities together is the first important step towards establishing more systematic thinking on ethics in mathematics, as explained further in \cite{Muller.2022}. 

Ethics in mathematics gives us a way to use our mathematics, and our mathematical insight, to make a tangible, positive impact on society, often by bringing our skills and insights to other areas, people and disciplines. Without a ``translator'' to explain what is happening in the mathematical community, the ability of external experts (e.g.~lawmakers regulating algorithms) is severely curtailed. Similarly, communicating external contributions to mathematicians is equally important. Mathematicians working on ethics in mathematics can act as ``bridges'' between those doing mathematical work, and those considering its effects on society. This work is very interesting in itself, bringing meaning and worth to our mathematical understanding, in addition to being an area of active recruitment in industry and beyond. 

Finally, there is the interesting challenge of teaching ethics in mathematics. Consider the cryptography problem posed in \cite{chiodo2023teaching}: Alice and Bob are using RSA encryption with private keys $(N, e_1)$ and $(N, e_2)$ with the same modulus $N$, and have asked their colleagues to send sensitive messages to each of them encrypted to their respective keys. An eavesdropper Eve is monitoring their communications. Using number theory we can show that Eve can always decrypt the messages. So who should we tell, and what should we tell them? What if it turns out that Alice and Bob have an arrest warrant out for them, and Eve is a police analyst; does that change our answer? What if the warrant is from a country with history of human rights violations; does the answer change now? Alice, Bob, and Eve may not know that decryption can occur, but as mathematicians we do. So how might we teach the consequences of our unique mathematical insight?

\newpage
\section*{Why is EiM tractable?}

Progress in mathematics is more than the number of theorems and published papers we produce; it is, as Thurston describes, the ``advance of human understanding of mathematics''\cite{WilliamP.Thurston.2}. Making progress on ethics in mathematics is also possible, and even goes beyond Thurston's criteria of ``advancing knowledge'' - it can be interpreted much more broadly by appreciating the theoretical and practical nature of the field. As a (relatively) nascent area of concern, ethics in mathematics has received comparatively little academic attention so far, thus there  remains much to be done. From the foundational level, where we need to get better at establishing and promoting the notion that there are indeed ethical issues that arise from mathematical work, to the practical applications of turning our theoretical knowledge into being morally good mathematicians and combatting morally bad uses of mathematics. Many important, critical problems need to be addressed, and therefore much work can, and should, be done.

We can't expect to ``solve'' ethics in mathematics with one publication, in the same way we can't solve topology with one publication. Our contributions are incremental, and, just like in mathematics, we make progress by standing on the shoulders of giants. There are many examples and case studies in this paper of  mathematics that raises ethical questions, many of which would benefit from further analysis, and there are even more that didn't fit into this paper or are yet to be found. There are rudimentary ethics guidelines and frameworks from various mathematical institutions and societies which are constantly being refined and updated (cf. \cite{Muller.2022, buell2022leveraging}). And there are introductory teaching modules on ethics in mathematics, and associated resources (e.g. \cite{chiodo2023teaching}), all of which could be substantially improved with additional time, effort, and insight. Work has already begun, but as mentioned in the previous section, it needs a larger community effort to make a more substantial impact.

Some of the problems associated with ethics in mathematics may seem daunting. The proliferation of AI through all aspects of society, and its many drawbacks and risks, is vast. The ever-increasing production of complicated financial instruments, and the threat they may pose to the global economy and to wealth inequality, might seem overwhelming. As mathematicians, we can't solve any of these problems by ourselves; they're simply too big. But by working on ethics in mathematics, we \textit{can} find ways to make a positive contribution, and thus make the world a little bit better. And in particular, we can contribute in a unique way and address issues only mathematicians are properly placed to address, or even identify! Ours might only be a small contribution, but it is critical, and one that only we can make.

In doing so, we can often draw ideas and inspiration from adjacent fields, including ``leveraging guidelines for ethical practice of [...] related professions''\cite[p.~5]{buell2022leveraging}. Engineering and computer science already have ways to recognise and deal with many of their ethical issues. Medicine has a long history of addressing ethics, and well-established ways to research and teach such problems. And some sub-disciplines of mathematics, such as statistics, have been addressing ethics for over a century. There are certainly lessons from there that, with a bit of work, can be carried over to mathematics more generally. As mathematics is deeply connected to the rest of the world, there are numerous open mathematical problems with strong connections to other existing and well-established ethics research areas. 

Ethics in mathematics has become an academic discipline, including peer-reviewed articles, books, conferences and dedicated courses\cite{ErnestHandbook}. Some large-scale funders of PhD studentships (e.g.~UK government Centres for Doctoral Training) now require dedicated ethical training to be incorporated in those centres, with other higher education agencies advocating for the inclusion of sustainability requirements into degrees due to the subject's increasing impact on society\cite{QAA2023}. Ethics in mathematics transcends academia, and is  critical to modern industries: the EU AI Act proposes new regulations for AI, and trustworthiness has become central in many areas, including finance, insurance, cryptography, and social media, all of which employ many mathematicians. 

Overall, ethics in mathematics is something we can do in a variety of ways, each giving clear avenues for us to make progress.

\newpage
\section*{Why is EiM non-trivial?}

As mathematicians we may be tempted to devalue ethics and think that, as a non-technical area, it is ``easy'' and ``obvious''. But neither is the case. Ethics, and in particular ethics in mathematics, is surprisingly challenging to work on. Often, it is only after seeing the ethical questions in one part of mathematics (e.g.~in our own area of expertise, or something receiving media coverage) that we start to grasp how ethics weaves itself through \textit{all} mathematics.

When dealing with ethics in mathematics, it's not only that the answers aren't obvious; even the very questions may not be evident. When setting up an optimisation problem, is it obvious that we should consider an objective function other than time or money? And what about the constraints; is it obvious that we should ask about aspects we were not initially briefed on? We know this from our mathematical research: being able to ask the right questions is often more than half the game. This is why mathematicians, and mathematics, are needed in the study of its ethics; in most cases, it can't be done without mathematical insight. We're the only ones who know the true significance of adding more constraints to a problem, so we're often the only ones who can ask about them. And as for deciding what the ``right'' goal to optimise for is, that's often a highly contentious, and thus non-trivial, question in itself. 

Operations research (OR) has begun developing ``soft OR techniques'' that help its practitioners to establish the ethics and responsibilities implicit in the OR process, with some even calling for a Hippocratic oath due to the potential for harm from bad models. But ethics in mathematics goes further, and at some point involves reverse engineering the process of carrying out mathematics to find out how to do it responsibly\cite{chiodo2023manifesto}. In doing so, the answers we find can take us completely by surprise. For example: a Hippocratic oath for mathematicians may be necessary, but such an oath would likely not be sufficient to address the shortcomings in the ethical awareness of ethically untrained mathematicians\cite{Muller.2022}. 

But it's not just the ethical problems themselves that pose a challenge; bringing together the right people to address these can also be challenging. Ethics in mathematics can't be done without help, support, and contribution from others in the mathematical community\cite{buell2022leveraging}. This requires mathematicians and institutions to come together in the right way. But growing such a community is hard, as not everyone will agree on what should be done, let alone how. Assembling people without due care can lead to conflict. But only assembling people who always agree with each other is no solution either, as we may end up lacking diversity and thus hinder ourselves from asking the right questions. Community building is one of the open challenges of doing ethics in mathematics effectively, with no straightforward answer.

Insight and communication are thus key challenges of ethics in mathematics. Even if we identify a problem in our, or someone else's, mathematical practice, we cannot just walk into the room and scream ``everything is on fire''. We need to offer pathways towards a solution, or we risk being ignored completely. Offering insight into what is, or has, gone wrong in the mathematical process (e.g.~inaccurate assumptions in a model) is a good start, but suggesting ways to resolve them is often more appreciated (e.g.~suggesting to consult domain experts to discuss how well our assumptions match reality). This requires astuteness to detect the problem, insight to realise what went wrong, foresight to see what needs to be done, and perspective to see how the proposal will be received by others, some of whom may not necessarily understand all the mathematics. 

That ethics in mathematics is hard should not be surprising. Many of its problems require a good grasp of many hard areas, including mathematics, ethics, philosophy and sociology of science, as well as people and context. All of these are far from trivial. We quickly realize that ethics in mathematics is neither one question nor one answer, but rather is woven through the entire mathematical process, so we cannot just ``do the ethics'' in one afternoon at the start of a project. 

Working on ethics in mathematics helps us appreciate that mathematics is everywhere and connects to all parts of life,  though not necessarily always for good. This is a challenge for us mathematicians in itself, as we may desire to see mathematics as morally pure, or neutral at worst\cite{Ernest2020}.

\newpage
\section*{Why is EiM well-defined?}
As early as primary school, we are taught that mathematics is certain and that its truth status cannot be influenced by external powers and interests\cite{BorbaSkovsmose1997}. In contrast, we encounter many debates and arguments about ethics which just appear to be exchanges about matters of opinion. However, a closer look at ethics and mathematics quickly reveals that the situation is more complex\cite{Nickel.2022}. Reasoning in ethics can be just as rational as reasoning in mathematics. It's just that its starting points are often less clear, and rarely communicated well in practice. 

Just like mathematics, ethical reasoning begins with finding good starting points or axioms, from which we can then build further. The starting points we find can be very different, depending on how we view and value different aspects of life. For some, this may mean being a utilitarian and calculating the good of an action by measuring its outcomes. Someone else may end up with a deontological position, applying universal rules to everyday ethics, while others may be more concerned with virtue, a duty of care for people, animals and the environment, or religious-centred ethics. The list of ethical positions is almost endless, and the positions often partially overlap but may also regularly disagree about what it means to be, and do, good. However, what they all have in common is that they produce well thought out answers subject to an initial set of assumptions; something we are all very familiar with as mathematicians. We also find these philosophical positions in ethics in mathematics. And so, ethics in mathematics is not a nebulous, fuzzy concept. Rather, it is well defined, modulo establishing what we care about and how we think about it.

But this means that to do ethics in mathematics effectively, we need to engage with others, understand their points of view, and obtain their input on the problem. It is often necessary to work towards consensus on a particular issue or challenge. We might seek consensus from colleagues on what assumptions to make, what to optimise over, how to interpret and convey our results to others, or even whether we should be carrying out the work or project in the first place. We seek consensus for decisions and choices which no one person can, or should, make on their own. Consensus is again a well-defined notion, which as mathematicians we use all the time in scenarios such as multi-refereed papers, PhD examination boards, and hiring panels. 

As for identifying whose responsibility it is to do ethics in mathematics, it quickly becomes evident that (almost) nobody else is doing it. Most mathematics departments around the world have no ethics committees to approve research studies, or ethics officers to deal with inner-departmental questions. Often the only solution is to do it ourselves, but even then we don't have to reinvent everything. By now there are well-defined research areas and many standard methods within ethics in mathematics. We have already alluded to adapting existing guidelines and frameworks from other areas. But there is plenty to choose from within mathematics, too. For example, building on the vast literature on mathematics education and social justice to produce resources, exercises, and assessments to be incorporated into higher-level mathematics courses. And there are also many well-known problems: AI safety, cryptography, finance, mathematics in warfare, and many more.

A feasible workflow of ethics in mathematics is often quite similar to that of mathematical research: 
\begin{enumerate}
    \item  Start by finding an area of concern and interest, and begin to understand your ethical position. Where are you right now? What are your starting points?
    \item \vspace{-4pt} Review the existing literature, including literature not published by mathematics. A lot of ethics in mathematics is published by other people, including philosophers, science and technology scholars, and domain experts. This is particularly true for mathematics education, where historically the community of research mathematicians can be quite distinct from those specialising in education.
    \item \vspace{-4pt} Engage with your peers and begin tackling the problem: ethics in mathematics is often a collaborative effort and you'll quickly find that those from outside mathematics mentioned in (2) are often very happy to work with mathematicians on an ethics paper. 
\end{enumerate}

\newpage
\section*{Why is EiM good?}
Is it not simply a tautology to say mathematics is ``good''? Well, yes, but we can go further than just saying that it is doing good for good's sake, and articulate exactly what this good is, and how it comes about. We all have different perspectives on what makes something ``good'', and doing ethics in mathematics helps us explore these and improve on them.

By doing ethics in mathematics, we do direct good for mathematics as a body of knowledge and its research practices. We can combat plagiarism and publication irregularities, thus preserving our corpus of knowledge. We can ensure that jobs are appointed fairly to those best placed to preserve and further mathematics. We can prevent mathematics from being used incorrectly, or misused, and thus enhance public trust in our work and mathematics more generally. We can give ourselves a better perspective on the problems we are solving, providing us with greater challenges and pushing our mathematical output to be even better. And we can ensure that mathematical work produces the best possible outcomes in the places and situations where it is used. 

By doing ethics in mathematics, we do direct good for the mathematical community and other communities. We can help turn mathematics into a field that is open and welcoming to all, regardless of background or circumstance. We can help ensure that mathematics is inclusive and diverse, thereby encouraging as many as possible to come and join the discipline. We can teach in ways that all students are able to learn and get the chance to appreciate as much mathematics as possible. We can teach and train upcoming mathematicians to be aware of the role, importance, and responsibility they have to other communities when carrying out their mathematical work. We can work towards ensuring that mathematics, as a tool and an output, does not discriminate or harm those individuals or groups who are most vulnerable in society. We can create mathematical tools that actively  combat existing discrimination and oppression in the world. And we can make progress to ensure that mathematics works for the benefit of all. Ethics in mathematics allows us to identify the most urgent problems and enables us to properly work towards a solution using our unique skillset.

By doing ethics in mathematics, we do direct good for the world at large. We can check to see what impact mathematics is having, and intervene in a preventative way before harm is caused. We can provide frameworks for mathematicians to follow, helping them to avoid producing mathematics which is detrimental to society. We can create tools to measure the performance and impact of mathematical output, to see where it might be causing harm in ways that are otherwise invisible to its creator. We can guide mathematicians to work on the right problems, and in the right ways. We can encourage mathematicians to factor in the socio-political context of their work,  to facilitate it being well-received rather than rejected. We can address the long-term effects of mathematics to society, the environment, and the planet. And we can promote humility in mathematicians, so that they better understand when they need to seek help and guidance from others, potentially outside  mathematics. Ethics in mathematics allows us to be more effective and efficient in our attempts to solve society's most pressing and difficult problems.

And by doing ethics in mathematics, we do direct good for ourselves as individual mathematicians. We can empower ourselves, allowing us to better navigate the world of mathematics and its culture. We expand our knowledge and expertise so that we can do more interesting and advanced work with better impact and minimal harm. We can extol the virtues of our work, making it more attractive for employers and funders. We can develop greater kindness and humility, towards those near to us, as well as those impacted far away whom we might never see. We can become part of new and different communities, interacting and working with people from other backgrounds and disciplines. We can find new opportunities for work and employment, bringing our now-broadened understanding of the role and place of mathematics into other domains and jobs. And we can be happy knowing that our work and input has made a positive contribution to the field, to those around us, to those in need, and to wider society.

Doing ethics in mathematics helps us to become a better mathematician: a better person doing better mathematics for a better world. 
\newpage
\section*{Going forward with EiM}

A useful observation, and subsequent resource, to have come out of the Ethics in Mathematics Project\footnote{Ethics in Mathematics Project: \url{https://www.ethics.maths.cam.ac.uk/}} at Cambridge University is the following decomposition of  ethical problems in mathematics, split up into what has been named as the ``10 pillars for responsible development''\cite{chiodo2023manifesto}: 

\vspace{-7pt}
\begin{enumerate}
    \setlength\itemsep{-3.5pt}
    \item Deciding whether to begin
    \item Diversity and perspectives
    \item Handling data and information
    \item Data manipulation and inference
    \item The mathematisation of the problem
    \item Communicating and documenting your work
    \item Falsifiability and feedback look
    \item Explainable and safe mathematics
    \item Mathematical artefacts have politics
    \item Emergency response strategies
\end{enumerate}

\begin{figure*}[hbt!]
    \centering
\tikzstyle{every node}=[draw=black, thick,anchor=west]
\tikzstyle{selected}=[draw=gray,fill=gray!30]
\begin{tikzpicture}[
  grow via three points={one child at (0.1,-0.7) and
  two children at (0.1,-0.7) and (0.1,-1.4)}, 
  edge from parent path={([xshift=10mm]\tikzparentnode.south west) |- (\tikzchildnode.west)}] 
  \node {Deciding whether to begin}
    child { node {Solve the right problem}}		
    child { node [selected] {Solve the problem the right way}
          child { node {How are you weighing up the consequences of your solution?}}
      child { node {Which issues are standard, obscure, or totally unexpected?}} 
      child { node [selected] {Why is a mathematical solution appropriate?}
            child{ node {How does mathematics actually help here?}}
            child{ node {Did you consider a non-mathematical approach?}}
            child{ node {Why is using mathematics the best option?}}
      }
    }
    child { node {$\dots$}};
\end{tikzpicture}
    \caption{The tree structure of ethics in mathematics (pillar 1)}
    \label{treeEiM}
\end{figure*}
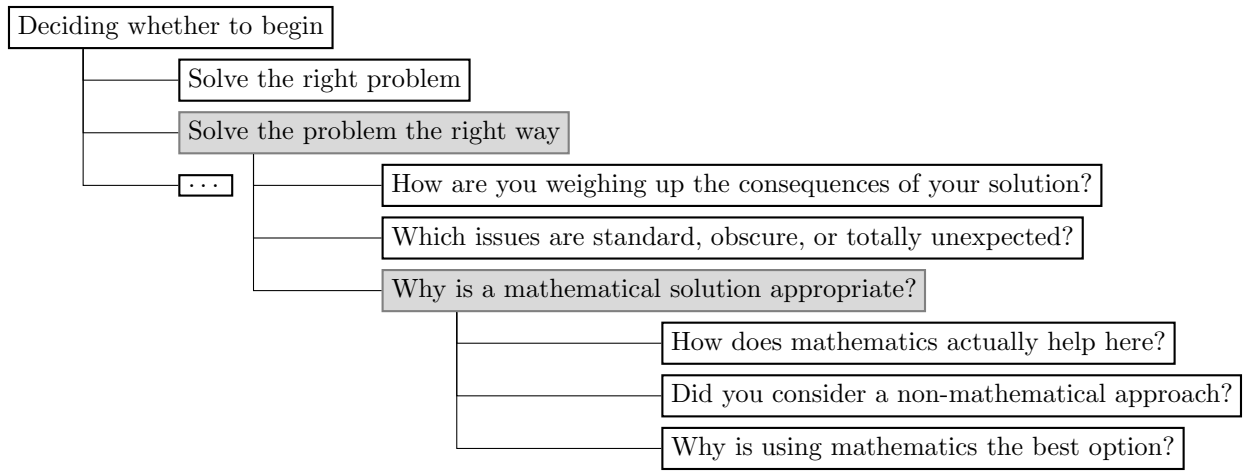

\vspace{-7pt}
Looking more closely, one observes in \cite{chiodo2023manifesto} that each ``pillar'' decomposes into 2--4 step-like tasks. Each task consists of 3--5 key questions, each of which breaks into a further 3--5 checkpoints, and each of those decomposes into a small set of actions (see figure \ref{treeEiM}). This demonstrates a clear tree-like, hyperbolic structure on how we might view ethics in mathematics; the more we explore a particular point or aspect, the more it breaks apart into several sub-aspects to consider, and this division may continue many times over. This can be seen as an \textit{opportunity} for those of us wishing to work on ethics in mathematics, as we can pick our desired level of specialisation and comfort, and there are open problems at every level of the tree.

Just like ethics in mathematics education\cite{DUBBS.2020}, the ethics of mathematical practice is done most effectively if it isn't done with a narrow definition of ethics in mind. Instead of being stumped by the question ``whose ethics am I supposed to consider?'', it is often better to acknowledge that the answer to this question may be found after working on an ethical question for a while. This is not surprising to us as mathematicians: only after studying an object and how it connects to other objects, do we establish how to properly define it and lay out its axioms. Writing down the hypotheses of a theorem is never the first step of its discovery, and neither is having an answer to the question ``whose ethics?'' a necessary requirement to start working on ethics in mathematics. Indeed, it's always part of exploring the issue at hand, and history has shown us that there can be more than one answer about the relationship between mathematics and ethics\cite{Nickel.2022}.

Ethics in mathematics - in research and practice - is a spectrum of actions and roles where \textit{everyone} can find something important to contribute to. As with any area, we can look at the full breadth of a problem, we can look at a particular (technical) aspect of a problem in great depth, and we can do many things in between. Some may work on overwhelmingly large problems, while others may want to work on something smaller and more manageable. This can be done in a complimentary, non-conflicting way that respects others. And quite often, there are hidden connections between the seemingly small problems and the biggest challenges. It's important not to excoriate someone's work on a problem just because we perceive a different problem as more urgent or important. We may not properly understand why the problem was chosen, or may not be aware of its depth and connections. Ethics in mathematics is often a lesson in kindness and humility. It requires from its participants an understanding that mathematics and its practice can change over time and across cultures and groups\cite{DAmbrosio1985}, and that its ethics may follow similar contingent circumstances. 

This work employs a vast array of methods and is inherently interdisciplinary and multidisciplinary; we need a lot of different types of understanding to be able to answer all the different types questions!  And any personal goal, and personal work preference, can fit into ethics in mathematics. We can do deep mathematics that aims to make the world better (e.g.~algorithmic fairness). We can do ethics that makes our mathematics better (e.g.~ensuring sufficient perspective in a team). We can apply elements from education, psychology, sociology and other social sciences or humanities to gain insight into the ethical issues of mathematical work (e.g.~to understand how existing ways of teaching affect our students' perception and use of mathematics). We can choose to treat ethics in mathematics as an oil-in-water suspension, where we add some mathematical insight into a predominantly ethical problem. Or we can treat it as water-in-oil, where we add ethical insight into a predominantly mathematical problem. Both of these are useful and productive contributions. Such work might not be 100\% mathematics, and may lack the ``exactness'' of doing standalone abstract mathematics, but that should not be a reason to disregard it. We often need to deal with non-mathematical factors, but these \textit{enhance} our mathematics by providing motivation for problems to work on, adding challenging constraints, and giving our output purpose and meaning. 

Starting this does not require any formal education in ethics, but simply for us to approach the subject and our work with an open mind and a willingness to do better. This is where looking at neighbouring subjects and debates can be very helpful. Many of us have some question or problem in  mind, but simply don't know how to start working on it. This can be in our teaching, in our department, in our research, or somewhere else entirely. Everyone starts doing ethics in mathematics with very little knowledge, but if we are commencing it today we have a clear advantage over those who worked on it 30, 10, or even 5 years ago. While Hersh deemed ``research work [to be] almost devoid of ethical content'' in his initial contributions to ethics in mathematics in the 1990s\cite[p.~23]{Hersh1990}, times have changed, and advances in knowledge have shown that there are ethical questions in \textit{all} areas of mathematics, both applied and pure. As a community, we are only just beginning to explore and understand these, but nobody has to start from zero anymore.

We may feel uncomfortable and tempted into inaction by thinking ``But I’m not trained to do this; I'm not an expert”. Ethics in mathematics is a very new area, and so there are few experts in the same way that one might be an expert in PDEs. This should not put us off, as we can still make a meaningful contribution despite being relatively inexperienced in the subject matter. Anyone spending a decent amount of time, probably just a few months, looking into ethics in mathematics would develop the necessary knowledge to make a novel and valuable contribution. We do not need to be the best in the world, nor have 20 years of experience, to work on ethics in mathematics. And that goes for teaching ethics in mathematics as well. We often teach undergraduate courses in areas where we are not internationally-recognised experts, but nonetheless know enough about, and sufficiently more than the students, to teach them. Ethics in mathematics need not be any different.

Many, though certainly not all, of the ethical concerns in mathematics are not specific to mathematics. They could just as easily apply to fields such as computer science and engineering, or even apply more generally to all vocations and lines of work. But that is not a reason to dismiss them. If such issues apply to other mathematical disciplines or even to \textit{everyone}, then they must surely apply to \textit{us} as well. And as outlined earlier, this means we can draw insight from other fields and carry those over in a more mathematics-centric way, to be maximally relevant for mathematicians.

Whatever our concerns are about mathematics, be they about continuing the body of knowledge, about caring for our mathematical and other communities, about the impact of mathematics on wider society, or indeed anything else, ethics in mathematics presents itself as a useful and good mechanism to highlight and address them. To us, the authors of this article, the study and practice of ethics in mathematics will only become more important in light of the ever-increasing use of mathematics. Only if we, as a community, reflect on and engage with the ethical questions that our field raises, can we ensure that its impact remains positive in the future.

\bibliography{Refs}

\end{document}